\documentclass{amsart}
\usepackage{amsfonts, amsthm, amssymb, amsmath, stmaryrd}
\usepackage{mathrsfs,array}
\usepackage{eucal,color,times,enumerate,accents}
\usepackage{tikz-cd}
\usepackage{bbm}
\usepackage[utf8]{inputenc}
\usepackage{pgfplots}
\pgfplotsset{compat=newest}
\usepackage{lipsum}
\usepackage{graphicx}
\usepackage{rotating}
\usepackage{enumitem}
\usepackage{enumerate}
\usepackage{footmisc}
\usepackage{textcomp}
\usepackage[bookmarksnumbered,bookmarksopen]{hyperref}
\hypersetup{colorlinks, linkcolor=black, citecolor=black}
\usepackage{cleveref}
\usepackage[colorinlistoftodos]{todonotes}
\usepackage[left=4cm,right=4cm,top=3cm,bottom=3cm]{geometry}
\usepackage{mathtools}
\newtheorem{theorem}[equation]{Theorem}

\newtheorem{lemma}[equation]{Lemma}

\newtheorem{corollary}[equation]{Corollary}

\numberwithin{equation}{section}

\theoremstyle{definition}

\newtheorem{question}[equation]{Question}

\newtheorem{remark}[equation]{Remark}

\usepackage{tikz}
\usepackage{tikz-cd}
\usepackage{rotating}
\usetikzlibrary{arrows,shapes,positioning}
\usetikzlibrary{decorations.markings}
\tikzstyle arrowstyle=[scale=1]
\usepackage{lipsum}

\usepackage{url}
\usepackage{imakeidx}
\makeindex[title=Index of definitions and notations,intoc]
\newcommand{\Z}{\mathbb{Z}}
\newcommand{\Q}{\mathbb{Q}}
\newcommand{\Oo}{\mathcal{O}}

\newcommand{\Ker}{\mathrm{Ker}}
\newcommand{\Image}{\mathrm{Im}}

\newcommand{\Spec}{\mathrm{Spec}}

\newcommand{\crys}{\mathrm{crys}}

\newcommand{\dr}{\mathrm{dR}}
\newcommand{\calH}{\mathcal{H}}
\newcommand{\fl}{\mathrm{flat}}
\newcommand{\et}{\mathrm{\acute et}}

\newcommand{\Kum}{\mathrm{Kum}}
\newcommand{\Cl}{\mathrm{cl}}

\newcommand{\Pp}{\mathbb{P}}

\newcommand{\C}{\mathbb{C}}
\newcommand{\Br}{\mathrm{Br}}
\newcommand{\NS}{\mathrm{NS}}
\newcommand{\dlog}{\mathrm{dlog}}

\title{Reduction modulo $p$ of the Noether's Problem}
\begin{document}
\date{\today}
\author{Emiliano Ambrosi and Domenico Valloni}
\begin{abstract}
Let $R$ be a complete valuation ring of mixed characteristic $(0,p)$ with algebraically closed fraction field $K$ and residue field $k$. Let $X/R$ be a smooth projective morphism. We show that if $X_k$ is stably rational, then $H^3_{\et}(X_K, \hat{\Z})$ is torsion-free. The proof uses integral $p$-adic Hodge theory of Bhatt-Morrow-Scholze and the study of differential forms in positive characteristic. We then apply this result to study the Noether problem for finite $p$-groups. 
\end{abstract}
\maketitle

\tableofcontents
\section{Introduction}
Many techniques to study the birational geometry of varieties are based on specialization methods in equicharacteristic or in mixed characteristic (see e.g. \cite{Voi15}, \cite{CP16}), often exploiting pathological phenomena in positive characteristic (see e.g. \cite{MR1273416, Tot16, MR4334404}). 

In this paper we study birational invariants in mixed characteristic. Let us fix a complete valuation ring $R$ of mixed characteristic $(0,p)$ with algebraically closed fraction field $K$ and residue field $k$ (necessarily algebraically closed) and a smooth proper morphism $X/R$. Motivated by the different behaviour of the Noether problem in positive and zero characteristic (see Section \ref{sec : appnoe} for more details) we study the following question:
\begin{question}
What can one say about the generic fiber $X_K$ knowing that the special fiber $X_k$ is stably rational (i.e. $X_k \times \Pp^n_k$ is birational to $\Pp^N_k$ for some $n,N\in \mathbb N$)?
\end{question} 
 In general, $X_K$ need not be stably rational, as we show for instance in Section \ref{sec : ex}, where we follow a suggestion of Colliot-Th\'{e}l\`{e}ne and adapt the constructions of \cite{MR3849287} to construct smooth proper morphisms $X/R$ with $X_K$ stably irrational and $X_k$ rational. On the other hand, we show that the following vanishing holds:
\begin{theorem} \label{main thm}
Let $X/R$ be a smooth proper morphism, and assume that the special fiber $X_k$ is stably rational. Then the $p$-torsion of $H^3_{\et}(X_K, \Z_p)$ vanishes.
\end{theorem}
If $K \subset \C$ is an embedding, then the vanishing of $H^3_{\et}(X_K, \Z_p)[p]$ together with smooth-proper base change in \'{e}tale cohomology (with $\ell$-adic coefficients, $\ell \neq p$) implies that $H^3_{\mathrm{sing}}(X(\C), \Z)$ is torsion free whenever $X_k$ is stably rational. We recall that the torsion of $H^3_{\mathrm{sing}}(X(\C), \Z)$ is a stably birational invariant of smooth proper varieties over $\C$, which is naturally isomorphic to $\Br(X)/\Br(X)_{\mathrm{div}}$, where $\Br(X):=H_{\et}^2(X,\mathbb G_m)$. This group was used by Artin and Mumford in their seminal paper \cite{MR321934} to give the first elementary example of a unirational threefold which was not rational. 
\subsection{Strategy}\label{sec : strat}
\numberwithin{equation}{subsection} 
To explain the strategy, we begin by showing an analogous statement concerning global differential forms. We retain the notation from the previous section. By the Kunneth formula and Hartogs' lemma, the vector spaces $H^0(X_K,\Omega^i_{X_K})$ are stably birational invariants of smooth proper $K$-varieties, and therefore they vanish if $X_K$ is stably rational. If $X/R$ is a smooth proper morphism, the semi-continuity theorem yields the inequality
\begin{equation} \label{eq : semi continuity}
    \dim_k(H^0(X_k,\Omega^i_{X_k}))\geq \dim_K(H^0(X_K,\Omega^i_{X_K})),
\end{equation}
from which it follows that if $X_k$ is stably rational, then $H^0(X_K,\Omega^i_{X_K})=0$ necessarily. 

Concerning \'{e}tale cohomology, as already mentioned, the proper smooth base change asserts that for a prime $\ell \neq p$ one has
$$H^3_{\et}(X_{K},\Z_{\ell})[\ell]\simeq H^3_{\et}(X_{k},\Z_{\ell})[\ell]=0,$$
where the last equality follows from the stable birational invariance of $H^3_{\et}(X_{k},\Z_{\ell})[\ell]$ and the fact that $X_k$ is stably rational. For $p$-adic coefficients, we can replace the smooth proper base-change with the integral $p$-adic Hodge theory of Bhatt-Morrow-Scholze, which will play the role of the semi-continuity theorem (\ref{eq : semi continuity}). By \cite[Theorem 1.1 (ii)]{BMS18}, one has the following inequality
\begin{equation}\label{eq : inequality}
\dim_k(H^3_{\crys}(X_{k}/W))[p])\geq \dim_{\mathbb F_p}(H^3_{\et}(X_{K},\mathbb Z_p)[p]),
\end{equation}
where $H^3_{\crys}(X_{k}/W)$ is the third crystalline cohomology group with integral coefficients of $X_{k}$. Thus, it would be enough to show that $\dim_k(H^3_{\crys}(X_k)[p])$ is a stably birational invariant of smooth proper varieties. It is unclear how to prove this without assuming resolution of singularities, also because crystalline cohomology behaves badly for open or singular varieties. We prove instead the following vanishing, which is enough to deduce Theorem \ref{main thm}:
\begin{theorem}\label{thm : maincrystalline}
Let $k$ be an algebraically closed field of characteristic $p$ and $X$ be a smooth proper $k$-variety.
Assume that 
\begin{enumerate}
	\item $H^i(X, \Oo_X) = 0$ for $ i= 2,3$;
 \item $H^0(X, \Omega^2_X) = 0$;
	\item $\Br(X)[p]=0$.
\end{enumerate}
Then $H^3_{\crys}(X/W)[p]=0$.
\end{theorem}
The assumptions of Theorem \ref{thm : maincrystalline} are satisfied if $X$ is stably rational, since all the conditions are stably birational invariants of smooth proper varieties (see \cite[Theorem 1]{CR11} for (1), point (2) follows from Hartog's lemma and  \cite[Corollary 6.2.11, Pag. 170]{MR4304038} gives (3)). Hence, Theorem \ref{main thm} follows from Theorem \ref{thm : maincrystalline} and (\ref{eq : inequality}).
\subsection{Applications to the Noether problem}\label{sec : appnoe}
\numberwithin{equation}{subsection}
As already mentioned, our original motivation is found in the Noether problem, which we now briefly recall. For a finite group $G$ and a field $K$, the Noether problem asks whether $\Pp^n_K/G$ is a stably rational variety, where $G$ acts on $\Pp^n_K$ in a linear and faithful way. The problem is well-posed since the stably birational class of the quotient does not depend on the chosen representation (\cite[Lemma 1.3]{Bog85}). For $K = \C$ we can then define the Artin-Mumford invariant $\mathrm{AM}(G)$ of $G$ as $\mathrm{Tors}(H^3_{\mathrm{sing}}(X,\mathbb \Z))$ where $X/ \C$ is any smooth proper birational model of $\Pp^n_K/G$. 
 
 The first counterexample of the Noether problem over $\C$ was given by Saltman in \cite{Salt84}, who for any prime $p$ produced a finite $p$-group $G$ for which the quotient $\Pp^n_{\mathbb C}/G$ is not stably rational, by showing that $\mathrm{AM}(G) \neq 0$. A general formula for $\mathrm{AM}(G)$ was later given by Bogomolov \cite[Theorem 3.1]{Bog88}, and this group is now known as the Bogomolov multiplier of $G$. 

The connection to Theorem \ref{main thm} comes from the classical observation that if $G$ is a $p$-group and $K$ has characteristic $p$, then $\mathbb P^n_K/G$ is always rational, see, e.g. \cite{Kun54} and \cite{zbMATH03141361}.

Keeping the notation as in the previous section, we can fix a finite $p$-group $G$, and linear faithful actions of $G$ on $\Pp^n_k$ and $\Pp^n_K$. Theorem \ref{main thm} together with the observation above implies the following: 
\begin{corollary}\label{cor : maincorollary}
If $\mathrm{AM}(G) \neq 0$, there does not exists a smooth proper $X/R$ such that $X_k$ is stably birational to $\Pp^n_k / G$ and $X_K$ is stably birational to $\Pp^n_K / G$.\end{corollary}
In particular, one sees that all the examples constructed in \cite{Salt84, Bog88} cannot have good stably rational reduction. This also implies that if $\mathrm{AM}(G) \neq 0$ and $G \rightarrow \mathrm{PGL}_n(R)$ is a representation such that the reduction map $G \rightarrow \mathrm{PGL}_n(k)$ is injective, then one cannot resolve the singularities of $\Pp^n_R / G$ relatively to $R$, i.e., there cannot be a smooth projective $X/R$ with a $R$-morphism $\pi \colon X \rightarrow \Pp^n_R / G$ which induces a resolution of singularities on both fibers. 

Peyre \cite{Pey08} constructed groups $G$ such that $\mathbb P^n_{\mathbb C}/G$ is not rational but $\mathrm{AM}(G)=0$. It is an interesting question at this point whether there is a $p$-group $G$ for which a resolution $X$ of $\Pp^n_R/G$ like in the corollary above exists, but $\Pp^n_K/G$ is not stably rational. 
\begin{remark}
A related phenomenon appears in the recent work of Lazda and Skorobogatov in \cite{zbMATH07694984}. They prove that, if $p=2$ and $Y\rightarrow \Spec(R)$ is an abelian surface such that $X_k$ is not supersingular, then one can resolve the singularities of $Y/ \{ \pm 1 \}$ to obtain a smooth proper morphism $X\rightarrow \Spec(R)$ such that the generic fiber is the Kummer variety $\Kum(Y_k)$ of $Y_k$ and the generic fiber is the Kummer variety $\Kum(Y_K)$ of $Y_K$. 

On the other hand, if $Y_k$ is surpersingular, then $\Kum(Y_k)$ is a rational surface due to \cite{Kat78}. Since $H^0(\Kum(Y_K),\Omega_{Y_K}^2)\neq 0$, the argument in the beginning of Section \ref{sec : strat} applies, and such $Y$ does not exist. Since $H^0(\Kum(Y_k),\Omega_{X_k}^2)=0$ if and only if $X_k$ is supersingular, the group $H^0(\Kum(Y_K),\Omega_{X_K}^2)$ is the only obstruction to the construction of such $Y$ in this case.
 \end{remark}
\subsection{Acknowledgement}
The authors are grateful to Stefan Schreieder for having introduced the second named author to the Noether problem, to Giuseppe Ancona for many interesting discussions and comments on a preliminary version of the article, and to Colliot-Th\'{e}l\`{e}ne for pointing out 
the example to us in Section \ref{sec : ex}. The authors would also like to thank both referees for their comments and for having greatly improved both the statement and the proof of our main result. 
\section{Proof of Theorem \texorpdfstring{\ref{thm : maincrystalline}}-}
In this section we prove Theorem  \ref{thm : maincrystalline}.
\subsection{Notation}
Let $k$ be an algebraically closed field of characteristic $p>0$ and let $X/k$ be a smooth proper variety such that
$$H^2(X,\mathcal O_X)=H^3(X,\mathcal O_X)=H^0(X,\Omega^2_{X/k})= \Br(X)[p]=0.$$
We let $\Omega_{X}^{\bullet}$ be the de Rham complex of $X$ and we define the following sheaves over $X$:
$$Z^i_X:=\Ker(d:\Omega_X^i\rightarrow\Omega^{i+1}_X);\quad B_X^i:=\Image(d:\Omega^{i-1}_X\rightarrow\Omega^{i}_X);\quad 
\calH^i_X=\frac{Z^i_X}{B_X^i}.$$

For every complex of sheaves $\mathcal F^{\bullet}$ over $X$ and every $i\in \mathbb N$, we let $\tau_{\geq i}\mathcal F^{\bullet}$ (resp. $\tau_{\leq i}\mathcal F^{\bullet}$)  be the upper (resp. lower) canonical truncation of $\mathcal F^{\bullet}$ and $\mathcal F^{\geq i}$ (resp. $\mathcal F^{\leq i}$) the upper (resp. lower) naive filtration of $\mathcal F^{\bullet}$. Recall that, for every $i\in \mathbb N$ there exists an exact triangle (see e.g.  \cite[Remark 08J5]{stacks-project}): 
$$\tau_{\leq i}\mathcal F^{\bullet}\rightarrow \mathcal F^{\bullet}\rightarrow \tau_{\geq i+1}\mathcal F^{\bullet}.$$
 \subsection{Preliminary reductions}
The universal coefficient theorem for crystalline cohomology (see e.g. \cite[(4.9.1), pag. 623]{illusiederhamwitt}), gives us an exact sequence
$$0\rightarrow H^2_{\crys}(X/W) \otimes k \rightarrow H^2_{\dr}(X)\rightarrow  H^3_{\crys}(X/W)[p]\rightarrow 0,$$ 
so that, to prove Theorem \ref{thm : maincrystalline}, it is enough to show that the natural map  
$$H^2_{\crys}(X/W) \otimes k \rightarrow H^2_{\dr}(X)$$ is surjective. 
In fact, we shall prove that the cycle class map $\Cl_{\dr} \colon \NS(X) \otimes k \rightarrow H^2_{\dr}(X/k)$ is surjective, which is enough due to the commutative diagram
\begin{equation}\label{eq : commutative cycles}
	\begin{tikzcd}
\NS(X)\otimes W\arrow{r}{\Cl_{\crys}}\arrow[two heads]{d} & H^2_{\crys}(X/W)\arrow{d}\\
\NS(X)\otimes k\arrow{r}{\Cl_{\dr}} & H^2_{\dr}(X),\\
	\end{tikzcd}
\end{equation}
where the first horizontal arrow is the crystalline cycle class map.
\subsection{Factorisation of the cycle class map}
To study $\Cl_{\dr}: \NS(X) \otimes k \rightarrow H^2_{\dr}(X/k)$, we factorise it in three arrows. Let $\iota: Z^1_X[-1]\rightarrow \Omega^\bullet$ be the natural inclusion. 

Recall that, by construction, the cycle class map  $\Cl_{\dr}: \NS(X) \otimes k \rightarrow H^2_{\dr}(X/k)$ is induced by the dlog map
$$\dlog: \mathcal O^*_X[-1]\rightarrow \Omega_X^\bullet$$
(see e.g. \cite[Section 0FLE]{stacks-project}), which factors trough
$$\mathcal O^*_X[-1]\xrightarrow{\dlog} Z^1_X[-1]\xrightarrow{\iota} \Omega_X^\bullet.$$

 Hence, the cycle class map factors trough the induced map $\iota:H^1(X,Z^1_X)\rightarrow H^2_{\dr}(X)$, giving a first factorisation 
$$\Cl_{\dr}:\NS(X) \otimes k\rightarrow  H^1(X,Z^1_X) \xrightarrow{\iota} H^2_{\dr}(X/k).$$ 

To go further, recall from \cite[(5.1.4) pag. 626]{illusiederhamwitt}, that for every $i$, there is a canonical isomorphism 
$H^i_{\fl}(X,\mu_p)\xrightarrow{\simeq} H^{i-1}(X,\mathcal O_X/(\mathcal O^*_X)^p)$ hence a canonical map
$\alpha: H^{2}(X,\mu_p)\rightarrow H^1(X,Z^1_X)$ induced again by $\dlog: \mathcal O_X/(\mathcal O^*_X)^p\rightarrow Z^1_X$. 
In conclusion, the map $\Cl_{\dr}:\NS(X) \otimes k\rightarrow  H^2_{\dr}(X)$ factorises further as 
$$\NS(X) \otimes k\xrightarrow{\Cl_{\fl}\otimes k} H^{2}_{\fl}(X,\mu_p)\otimes k\xrightarrow{\alpha} H^1(X,Z^1_X)\xrightarrow{\iota}H^2_{\dr}(X),$$
where $\Cl_{\fl}: \NS(X)\rightarrow H^{2}_{\fl}(X,\mu_p)$ is the cycle class map in flat cohomology. 
\subsection{Studying the factorisation}
To show that $\Cl_{\dr}:\NS(X) \otimes k\rightarrow  H^2_{\dr}(X)$ is surjective, we will show that:
\begin{enumerate}
	\item the map $\Cl_{\fl}:\NS(X) \rightarrow H^{2}_{\fl}(X,\mu_p)$ is surjective;
	\item the map $\alpha: H^{2}_{\fl}(X,\mu_p) \otimes k \rightarrow H^1(X,Z^1_X)$ is an isomorphism;
	\item the map $\iota: H^1(X,Z^1_X)\rightarrow H^2_{\dr}(X)$ is surjective.
\end{enumerate}
\subsubsection{Proof of (1)}
The map $\Cl_{\fl}:\NS(X) \rightarrow H^{2}_{\fl}(X,\mu_p)$ is induced by the connecting map 
$$H^1(X,\mathcal O_X^*)= \mathrm{Pic}(X)\rightarrow H^2_{\fl}(X,\mu_p)$$
in the long exact sequence associated to the short exact sequence $$0\rightarrow \mu_p\rightarrow \mathcal O_X^*\xrightarrow{(-)^p} \mathcal O_X^*\rightarrow 0$$
of sheaves in flat site of $X$. 
Hence the surjectivity of $\Cl_{\fl}:\NS(X) \rightarrow H^{2}_{\fl}(X,\mu_p)$ follows directly from the assumption $H_{\fl}^2(X,\mathcal O_X^*)=H_{\et}^2(X,\mathcal O_X^*)=Br(X)=0$. 
\subsubsection{Proof of (2)}
Recall from \cite[(2.1.23), pag. 518]{illusiederhamwitt} the exact sequence of \'etale sheaves
$$0\rightarrow \mathcal O_X^*/(\mathcal O_X^*)^p\xrightarrow{\dlog} Z^1_X\xrightarrow{i-C} \Omega^1_X\rightarrow 0,$$
where $i:Z^1_X\rightarrow \Omega^1_X$ is the natural inclusion and $C:Z^1_X\rightarrow \Omega^1_X$ is the Cartier operator. 
Since \'etale cohomology and Zariski cohomology agree for coherent sheaves,  taking the associated long exact sequence in cohomology, we get an exact sequence
\begin{equation*} \label{long sequence fppf}
H^0(X, Z^1_X) \xrightarrow{i-C} H^0(X,  \Omega^1_X) \rightarrow H^2_{\mathrm{fppf}}(X, \mu_p)  \rightarrow H^1(X, Z^1_X) \xrightarrow{i-C} H^1(X, \Omega^1_X)
\end{equation*}
where  $i$ is a linear morphism and $C$ is a Frobenius-linear map. 

To prove (2), it is then enough to show that 
$H^0(X, Z^1_X) \xrightarrow{i-C} H^0(X,  \Omega^1_X)$ is surjective and that 
$$\Ker(H^1(X, Z^1_X) \xrightarrow{i-C} H^1(X, \Omega^1_X))\otimes k\simeq H^1(X, Z^1_X).$$

Since $H^i(X, Z^1_X)$ and $H^i(X, \Omega^1_X)$ are finite-dimensional vector spaces, we can use the following classical lemma (see e.g. \cite[Lemma 4.13, pag 128]{milneEtale}).
\numberwithin{equation}{subsubsection}
\begin{lemma} \label{lemma semilinear}
Let $k$ be an algebraically closed field of characteristic $p>0$ and let $V$ be a finite-dimensional vector space over $k$. Let $f:V \rightarrow V$ be a $k$-linear isomorphism and let $C \colon V \rightarrow V$ be any Frobenius-linear map. Then:
\begin{enumerate}
	\item[(a)] $f - C$ is surjective
	\item[(b)] if in addition $C$ is bijective, then $V = \ker(f - C) \otimes k$
\end{enumerate}
\end{lemma}
\numberwithin{equation}{subsection}
Hence to prove (2) it is enough to show that 
$$i:H^0(X, Z^1_X) \rightarrow H^0(X, \Omega^1_X), \quad  i: H^1(X, Z^1_X) \rightarrow H^1(X, \Omega^1_X)\quad \text{and} \quad C: H^1(X, Z^1_X) \rightarrow H^1(X, \Omega^1_X)$$
are isomorphisms. 

The fact that the map $H^0(X, Z^1_X) \rightarrow H^0(X, \Omega^1_X)$ is an isomorphism, follows directly from the exact sequence of sheaves 
$$0\rightarrow Z^1_X\rightarrow \Omega^1_X\xrightarrow{d} \Omega^2_X$$
and the assumption $H^0(X, \Omega_X^2)=0$.

Hence, we are left to prove that 
$$i: H^1(X, Z^1_X) \rightarrow H^1(X, \Omega^1_X)\quad \text{and}\quad C: H^1(X, Z^1_X) \rightarrow H^1(X, \Omega^1_X) $$ are isomorphisms. By comparing dimensions, it is enough to prove:
\begin{enumerate}
	\item[(i)] $i: H^1(X, Z^1_X) \rightarrow   H^1(X, \Omega^1_X)$ is injective;
	\item[(ii)] $C: H^1(X, Z^1_X) \rightarrow   H^1(X, \Omega^1_X)$ is surjective.
\end{enumerate}
\proof[Proof of (i)]
Since $Z_X^1[-1] = \tau_{\leq 1} \Omega_X^{\geq 1}$, there is an exact triangle
$$Z_X^1[-1] \rightarrow \Omega_X^{\geq 1} \rightarrow \tau_{\geq 2}\Omega_X^{\geq 1},$$
so that the map $H^1(X, Z^1_X) \rightarrow H^2(X,\Omega_X^{\geq 1})$ is injective, since 
$H^1(X,\tau_{\geq 2}\Omega_X^{\geq 1})=0$, because 
$\tau_{\geq 2}\Omega_X^{\geq 1}$ is concentrated in degrees $\geq 2$. 
So it is enough to show that the natural map $H^2(X,\Omega_X^{\geq 1}) \rightarrow H^2(X,\Omega_X^1[-1])=H^1(X,\Omega_X^1)$ induced by the map $\Omega_X^{\geq 1} \rightarrow \Omega_X^1[-1]$ is injective. But the latter fits in the short exact sequence of complexes 
$$0 \rightarrow \Omega_X^{\geq 2} \rightarrow \Omega_X^{\geq 1} \rightarrow \Omega_X^1[-1] \rightarrow 0, $$
so that we just need to prove that $H^2(X, \Omega_X^{\geq 2})=0$. On the other hand, 
$$H^2(X, \Omega_X^{\geq 2})=\Ker(H^0(X, \Omega_X^{2}) \xrightarrow{d} H^0(X, \Omega_X^{3}))\subseteq H^0(X, \Omega_X^{2})=0$$
by our assumption on $X$. Hence  $i:H^1(X, Z^1_X) \rightarrow   H^1(X, \Omega^1_X)$ is injective.

	\proof[Proof of (ii)]: Since $k$ is perfect, from \cite[(2.1.22)]{illusiederhamwitt} we have a short exact sequence of sheaves
$$0\rightarrow B^1_X\rightarrow Z^1_X\xrightarrow{C} \Omega^1_{X}\rightarrow 0,$$
where $C$ is the Cartier operator.
Hence it is enough to show that $H^2(X, B^1_X)$ is zero. But this follows from the short exact sequence $$ 0 \rightarrow \Oo_X \xrightarrow{(-)^p} \Oo_X \xrightarrow{d} B^1_X \rightarrow 0 $$
and the assumption $H^2(X, \Oo_X) = H^3(X, \Oo_X) = 0$.
 \subsubsection{Proof of (3)}
 Since $H^2(X,\mathcal O_X)=0$ by assumption, the short exact sequence
 $$0\rightarrow  Z_X^1[-1]\rightarrow  \tau_{\leq 1}  \Omega_{X}^{\bullet}\rightarrow \mathcal O_X\rightarrow 0$$
 shows that the natural map $H^1(X, Z^1_X)\rightarrow  H^2(X,\tau_{\leq 1}  \Omega_{X}^{\bullet})$ is surjective. So, it is enough to show that the natural map 
 $H^2(X,\tau_{\leq 1}  \Omega_{X}^{\bullet}) \rightarrow  H^2(X,\Omega_{X}^{\bullet})=H^2_{\dr}(X)$ is surjective. Since there is an exact triangle
\begin{equation*}
   \tau_{\leq 1}  \Omega_{X}^{\bullet} \rightarrow  \Omega_{X}^{\bullet} \rightarrow \tau_{\geq 2} \Omega_{X}^{\bullet},
\end{equation*} 
it is then enough to show that $H^2(X, \tau_{\geq 2} \Omega_{X}^{\bullet})=0$. But 
$$H^2(X, \tau_{\geq 2} \Omega_{X}^{\bullet}) = \ker(H^0(X, \tau_{\geq 2} \Omega_{X}^{\bullet}[2])\rightarrow  H^0(X,\Omega^3_X))=H^0(X, \ker( \Omega^2_X / B^2_X \rightarrow  \Omega^3_X)) = H^0(X, \mathcal{H}^2_X).$$ 
Again, since $k$ is perfect, the inverse Cartier operator (\cite[(2.1.22)]{illusiederhamwitt}) gives an isomorphism $\mathcal{H}^2( \Omega^\bullet_X) \cong \Omega_{X}^2$, so 
$$H^2(X, \tau_{\geq 2} \Omega_{X}^{\bullet}) \simeq H^0(X,\Omega_{X}^2)= 0.$$ Hence the natural map $H^1(X, Z^1_X) \rightarrow H^2_{\dr}(X/k)$ is surjective.  This concludes the proof of (3) and the proof of Theorem \ref{thm : maincrystalline}.
\section{A stably irrational variety reducing to a rational variety} \label{sec : ex}
Let $R$ be the ring of integers of $K:=\C_p$ and $k$ its residue field. 
In this last section, we show how to construct, for every $p\gg 0$, examples of smooth proper schemes $X\rightarrow \Spec(R)$ such that $X_K$ is not stably rational and such that $X_k$ is rational, as suggested to us by Colliot-Th\'{e}l\`{e}ne.
The construction uses and is based on the analogous construction in \cite{MR3849287} of a family of proper smooth varieties over the complex number with stably irrational general fiber but with some rational fiber. 
\subsection{A general lemma}
We begin with a general lemma which reduces the construction of examples to the construction of mixed characteristic families with properties that are easier to check. 
Let $B/R$ be smooth with geometrically integral fibers and $X \rightarrow B$ a smooth proper family of varieties.
\begin{lemma}\label{lem : reduction}
Assume that there exists a point $b\in B(\C_p)$ such that $X_b$ is not stably rational. Then, for every $a\in B(k)$ there exists a lift $b'\in B(R)$ of $a$ such that $X_{b'}$ is not stably rational.  
\end{lemma}
\proof
 By \cite[Corollary 4.1.2]{NicaiseOttem}, the set 
$$ B(\C_p)_{r} \coloneqq \{ b \in B(\C_p) \colon X_b \text{ is stably rational} \}$$
is a countable union of closed subvarieties. Define now $B(\C_p)_{nr} \coloneqq B(\C_p) \setminus B(\C_p)_{r}$. By the assumption on $b$, the set $B(\C_p)_{r}$ is the countable union of \textit{proper} closed subvarieties. 

Since, by the Hensel lemma, the map $\pi:B(R) \rightarrow B(k)$ is surjective, we can choose a lift $b''$ of $a$. The set $\pi^{-1}(a)\subseteq B(R) \subseteq B(K)$ is an open neighborhood of $b''$ in $B(K)$.  Since $B(\C_p)_{r}$ is the countable union of proper closed subvarieties, we can apply \cite[Lemma 4.29]{zbMATH06073742} to deduce that there exists a $b'\in B(\C_p)_{nr}\cap \pi^{-1}(a)$. This concludes the proof.
\endproof
\subsection{An example}
By Lemma \ref{lem : reduction}, to construct a smooth proper scheme $X\rightarrow \Spec(R)$ such that $X_K$ is not stably rational and such that $X_k$ is rational, it is enough to construct a family $X\rightarrow B$ over $R$ such that there exist points $b\in B(\mathbb C_p)$ and $a\in B(k)$ such that $X_b$ is stably irrational and $X_a$ is rational.  Such a family can be constructed using directly Hassett, Pirutka, and Tschinkel example \cite{MR3849287}, see also \cite[Section 12.2.2]{MR4304038} and \cite{zbMATH07036862}. We give some details. 

Let $X'\rightarrow Z$ be the universal family of quadric bundles over $\mathbb P^{2}_{\Q}$ given in $\mathbb P_\Q^2\times \mathbb P_\Q^3$ by a bihomogeneous form of bidegree $(2, 2)$. After choosing coordinates $x,y,z$ and $U,V,W,T$ on $\mathbb P_\Q^2\times \mathbb P_\Q^3$, the variety $Z$ identifies with the space of bihomogeneous forms $F = F(x,y,z,U,V,W,T)$ of bidegree $(2,2)$ in $\Pp^2_\Q \times \Pp^3_\Q$ that are symmetric quadratic forms in the variables $U,V,W,T$, since any such $F$ determines a quadric bundle over $\Pp^2_\Q$ via the projection $\Pp^2_\Q \times \Pp^3_\Q \rightarrow \Pp^2_\Q$. In turn, these forms are given by a $4 \times 4$ symmetric matrix $A = (a_{i,j})_{1 \leq i,j \leq 4}$ where each entry $a_{i,j} = a_{i,j}(x,y,z)$ is a homogeneous polynomial of degree $2$.

By the arguments in  \cite{zbMATH07036862} and  Bertini theorem, there exists a dense open Zariski $B_{\Q} \subset Z$ such that the restriction of the family $X_{\Q}\rightarrow B_{\Q}$ to $B_{\Q}$ parametrizes quadric bundle flat over $\Pp^2_\Q$ and with smooth total space. 

By spreading out, this construction extends to give a smooth family $X \rightarrow B$ over $\mathbb{Z}[1/n]$ for $n$ large enough whose base change to $\Q$ identifies with $X_{\Q}\rightarrow B_{\Q}$. By the main result of \cite{MR3849287}, the general fiber of $X(\mathbb C)\rightarrow B(\mathbb C)$ is not stably rational, hence, for every $p$, there exists a $b\in B(\mathbb C_p)$ such that $Y_b$ is not stably rational. So, we are left to show that for $p\gg 0$, there exists a $a\in B(\overline {\mathbb F}_{p})$ such that $X_a$ is rational.

Using Bertini, there exists a rational point $r\in B(\overline{\Q})$ such that the corresponding $4 \times 4$ symmetric matrix $A = (a_{i,j})_{1 \leq i,j \leq 4}$ has $a_{1,1}=0$. By spreading out, we can choose $p\gg 0$ such that $r \in B(\overline{\Q})$ extends to a point $\widetilde a \in B(\overline{\Z}_p)$ whose reduction $a$ modulo $p$ defines a flat quadric bundle $X_a\rightarrow  \mathbb P^2_{\overline{\mathbb F}_p}$ with smooth total space, whose associated matrix has $a_{1,1}=0$.
Since $a_{1,1}=0$, for every $x\in \mathbb P^2_{\overline{\mathbb F}_p}$, the point $[1:0:0:0]$ is $k(x)$-rational point of the fiber of $X_a\rightarrow \mathbb P^2_{\overline{\mathbb F}_p}$ in $x$. In particular the morphism $X_a\rightarrow \mathbb P^2_{\overline{\mathbb F}_p}$ has a rational section, hence $X_a$ is rational. This concludes the construction of a proper smooth scheme over $R$ with rational special fiber and stably irrational generic fiber. 
\bibliographystyle{alpha}
\bibliography{bibliographyT.bib}
\vspace{5pt}
Emiliano Ambrosi, \textsc{University of Strasbourg}\par\nopagebreak
 \textit{E-mail address}: \texttt{ eambrosi@unistra.fr}\\
\vspace{1pt}
 Domenico Valloni, \textsc{University of Hannover}\par\nopagebreak
 \textit{E-mail address}: \texttt{ 
valloni@math.uni-hannover.de}

\end{document}